\documentclass[final]{elsarticle}%
\usepackage{amsfonts}
\usepackage{amsmath}
\usepackage{amssymb}
\usepackage{graphicx}%
\providecommand{\U}[1]{\protect\rule{.1in}{.1in}}
\journal{ArXiv}
\newtheorem{theorem}{Theorem}

\newtheorem{definition}[theorem]{Definition}
\newtheorem{example}[theorem]{Example}

\newtheorem{lemma}[theorem]{Lemma}

\newtheorem{proposition}[theorem]{Proposition}

\newenvironment{pf}[1][Proof]{\textbf{#1.} }{\  \rule{0.5em}{0.5em}}
\begin{document}
%
\begin{frontmatter}%


%

\title{Kneading determinants of infinite order linear recurrences}%

%

\author{João F. Alves\corref{cor1}\fnref{label1}}%
%

\ead{jalves@math.tecnico.ulisboa.pt}%
%

\cortext[cor1]{Corresponding Author}%

\fntext[label1]{Centro de Análise Matemática Geometria e Sistemas Dinâmicos, Math. Dep., Tech. Institute of Lisbon, Univ. of Lisbon, Av. Rovisco Pais, 1049-001 Lisbon, Portugal}%

\author{Ant\'onio Bravo\fnref{label2}}%
%

\ead{abravo@math.tecnico.ulisboa.pt}%

\fntext[label2]{Centro de Análise Funcional e Aplicações, Math. Dep., Tech. Institute of Lisbon, Univ. of Lisbon, Av. Rovisco Pais, 1049-001 Lisbon, Portugal}%

\author{Henrique M. Oliveira\fnref{label1}}%
%

\ead{holiv@math.tecnico.ulisboa.pt}%
%

\begin{abstract}%
Infinite order linear recurrences are studied via kneading matrices and
kneading determinants. The concepts of kneading matrix and kneading
determinant of an infinite order linear recurrence, introduced in this work,
are defined in a purely linear algebraic context. These concepts extend the
classical notions of Frobenius companion matrix to infinite order linear
recurrences and to the associated discriminant of finite order linear
recurrences. Asymptotic Binet formulas are deduced for general classes of
infinite order linear recurrences as a consequence of the analytical
properties of the generating functions obtained for the solutions of these
infinite order linear recurrences.%

\end{abstract}%
%

%

\end{frontmatter}%



\section{Introduction}

The concept of kneading determinant was introduced by Milnor and Thurston
\cite{knead} in the late eighties of the last century in the context of
one-dimensional dynamics. Later on, it was shown that the kneading determinant
of an interval map can be regarded as the determinant of a pair of linear
endomorphisms with finite rank, see \cite{jalves1} and \cite{jalves2}. This
latter point of view, purely linear algebraic, is the link between the Milnor
and Thurston notion and our definition of kneading determinant of a linear
recurrence. Indeed, as we will see, the kneading determinant of an infinite
order linear recurrence is a particular case of the above mentioned
determinant of a pair of linear endomorphisms with finite rank described in
\cite{jalves2}.

The main objective of this paper is to show that the kneading determinants
play an important role in the study of infinite vector recurrences, giving
directly the generating functions of the solution of the problem. In addition,
the determinants present a powerful computational tool to obtain the actual
solutions of finite and infinite order linear recurrences..

Linear recurrences have a long history, they constitute generalizations of the
eight centuries old finite linear recurrences of Leonardo de Pisa, or
Fibonacci \cite{fibo}%
\[
q_{n+1}=q_{n}+q_{n-1}\text{, with }q_{0}=0\text{, }q_{1}=1\text{ and }n\geq1.
\]

In the 19th century Jacques Philippe Marie Binet popularized a formula,
earlier known to De Moivre, solving the Fibonacci recurrence as a function of
$n$.

In a series of papers, \cite{rash,rash1,rashmotta2} Rachidi and other authors
studied linear infinite order scalar recurrences. Given an infinite sequence
of coefficients $\left\{  a_{i}\right\}  _{i=0,1,2\ldots}$, with some possible
conditions on the sequence, like periodicity \cite{rashmotta2}, positivity of
the coefficients, see the recent work \cite{rash2},\ or the existence of some
limit, the problem was to find a solution of the infinite order linear scalar
recurrences%
\begin{equation}
q_{n+1}=%
{\displaystyle\sum\limits_{i=0}^{+\infty}}
a_{i}q_{n-i},\text{ for }n\geq0\text{,}\label{Rashidi}%
\end{equation}
with an infinite set of initial conditions $\left\{  q_{i}\right\}
_{i=0,-1,-2,\ldots}$. By studying the results of these researchers, namely on
Binet formulas, we adopted a new approach to the problem using the different
technique of kneading determinants. We apply this new method to a wider class
of recurrences, obtaining solutions and asymptotic behaviour showing the
conceptual and computational power of kneading determinants. One of the
advantages of using generating functions is the possibility of analyzing the
asymptotic behaviour using the analytical properties of the generating function.

The paper is organized as follows, in sections $2$ and $3$, we introduce the
terminology and the main results of this paper, we generalize this problem to
vectorial recurrences (\ref{vector-rec}) and present their solutions.
Naturally, the method solves scalar recurrences as a particular case. We
present three fundamental results characterizing the solutions of infinite
order linear recurrences, Theorem \ref{T1} generalizes the concept of
Frobenius companion matrix, Theorem \ref{T2} gives the generating function for
the solutions of the recurrences and Theorem \ref{T3} gives asymptotic Binet
formulas for the asymptotic solutions of the problem. The technical details
and proofs are given in the last two sections.

\section{Terminology and definitions}

Let $\mathbb{N}=\left\{  0,1,2,...\right\}  $ be the set of non-negative
integers, $p$ a positive integer and $\left\{  \mathbf{A}_{n}\right\}
_{n\in\mathbb{N}}$ an infinite sequence of $p\times p$ matrices with complex
entries. In this paper we are interested in vectorial (or matricial)
homogeneous linear recurrences of the type%

\begin{equation}
\mathrm{x}_{n+1}=\sum\limits_{i=0}^{+\infty}\mathbf{A}_{i}\mathrm{x}%
_{n-i}\text{, for }n\in\mathbb{N}\text{,} \label{vector-rec}%
\end{equation}
where $\mathrm{x}_{n}\in\mathbb{C}^{p}$ for all $n\in\mathbb{Z}$ and
$\mathrm{x}_{n}=\mathrm{0}$ for almost all\footnote{In the sense of all but
except finitely many.} negative integers $n$.

We call this type of homogeneous linear recurrences of generalized Fibonacci
type on $\mathbb{C}^{p}$, for short, $Fib_{p}$ recurrences.

A $Fib_{p}$ recurrence is completely determined by a sequence of matrices
$\left\{  \mathbf{A}_{n}\right\}  _{n\in\mathbb{N}}$. If there exists $k\in$
$\mathbb{Z}^{+}$ such that $\mathbf{A}_{k-1}\neq\mathbf{0}$ and $\mathbf{A}%
_{n}=\mathbf{0}$ for $n\geq k$, the $Fib_{p}$ recurrence is said to be of
finite order $k$ (for short a $Fib_{p}^{k}$ recurrence). If the set $\left\{
n\in\mathbb{N}:\mathbf{A}_{n}\neq\mathbf{0}\right\}  $ is infinite, the
$Fib_{p}$ recurrence is said to be of infinite order (for short a
$Fib_{p}^{\infty}$ recurrence). With this notation the recurrence
(\ref{Rashidi}) is of type $Fib_{1}^{\infty}$ and the original Fibonacci
recurrence is of type $Fib_{1}^{2}$.

The concepts of kneading matrix and kneading determinant of a $Fib_{p}$
recurrence, introduced in this work, will be defined in a purely linear
algebraic context. These concepts extend for $Fib_{p}^{\infty}$ recurrences
the classical notions of Frobenius companion matrix and associated
discriminant of a $Fib_{p}^{k}$ recurrence.

Throughout the paper $\mathbb{C}[z]$ and $\mathbb{C[}[z]]$ denote respectively
the commutative rings of polynomials and formal power series with complex
coefficients. Matrices with entries in $\mathbb{C}$, $\mathbb{C}[z]$ and
$\mathbb{C[}[z]]$ will be denoted respectively as elements of $\mathbb{C}%
^{m\times n}$, $\mathbb{C}[z]^{m\times n}$ and $\mathbb{C}[[z]]^{m\times n}$.
The $m\times m$ identity matrix $\mathbf{I}_{m}$ will be usually written
$\mathbf{I}$ keeping in mind that its order is always well defined from the context.

The infinite-dimensional vector spaces over $\mathbb{C}$%
\begin{equation}
\boldsymbol{U}=\bigoplus\limits_{n\in\mathbb{N}}\mathbb{C}^{p}\text{ and
}\boldsymbol{V}=%
{\textstyle\prod\limits_{n\in\mathbb{Z}}}
\mathbb{C}^{p} \label{spaces}%
\end{equation}
will play an important role in this discussion. We write $\boldsymbol{u}$ and
$\boldsymbol{v}$ for denoting the vectors of $\boldsymbol{U}$ and
$\boldsymbol{V}$ with components $\mathrm{u}_{n}\in\mathbb{C}^{p}$ and
$\mathrm{v}_{n}\in\mathbb{C}^{p}$, i.e.,%
\[
\boldsymbol{u}=\left(  \mathrm{u}_{n}\right)  _{n\in\mathbb{N}}=(\mathrm{u}%
_{0},\mathrm{u}_{1},...)\text{ with }\mathrm{u}_{n}\in\mathbb{C}^{p}\text{ }%
\]
and%
\[
\boldsymbol{v}=\left(  \mathrm{v}_{n}\right)  _{n\in\mathbb{Z}}%
=(...,\mathrm{v}_{-1},\mathrm{v}_{0},\mathrm{v}_{1},...)\text{ with
}\mathrm{v}_{n}\in\mathbb{C}^{p}.
\]

In contrast with $\boldsymbol{V}$, the space $\boldsymbol{U}$ admits a
countable infinite basis. From now on we reserve the symbols $\mathrm{e}%
_{1},...,\mathrm{e}_{p}$ for denoting the vectors of the standard basis of
$\mathbb{C}^{p}$ and $\boldsymbol{e}_{\beta}$, with $\beta\in\mathbb{Z}^{+}$,
for denoting the vectors of the standard basis of $\boldsymbol{U}$:%
\[%
\begin{array}
[c]{l}%
\boldsymbol{e}_{1}=(\mathrm{e}_{1},\mathrm{0},\mathrm{0},.....),\boldsymbol{e}%
_{2}=(\mathrm{e}_{2},\mathrm{0},\mathrm{0},.....),.....,\boldsymbol{e}%
_{p}=(\mathrm{e}_{p},\mathrm{0},\mathrm{0},.....),\\
\boldsymbol{e}_{p+1}=(\mathrm{0},\mathrm{e}_{1},\mathrm{0},...),\boldsymbol{e}%
_{p+2}=(\mathrm{0},\mathrm{e}_{2},\mathrm{0},...),...,\boldsymbol{e}%
_{2p}=(\mathrm{0},\mathrm{e}_{p},\mathrm{0},...)\text{,}\\
\boldsymbol{e}_{2p+1}=(\mathrm{0},\mathrm{0},\mathrm{e}_{1},\mathrm{0}%
,...),...
\end{array}
\]
where $\mathrm{0}$ denotes the zero vector of $\mathbb{C}^{p}$.

After these basic remarks we now present the terminology of linear recurrences.

A vector $\boldsymbol{v}=\left(  \mathrm{v}_{n}\right)  _{n\in\mathbb{Z}}$
$\in\boldsymbol{V}$ is said to be a solution of a $Fib_{p}$ recurrence
(\ref{vector-rec}) if the set $\left\{  n<0:\mathrm{v}_{n}\neq\mathrm{0}%
\right\}  $ is finite and
\[
\mathrm{v}_{n+1}=\sum\limits_{i=0}^{+\infty}\mathbf{A}_{i}\mathrm{v}%
_{n-i},\text{ for all }n\geq0.
\]
The subspace of $\boldsymbol{V}$ whose vectors are the solutions of the
$Fib_{p}$ recurrence is denoted by $\boldsymbol{S}$.

Naturally, there exists an isomorphism%
\[%
\begin{array}
[c]{cccc}%
\Theta: & \boldsymbol{U} & \rightarrow & \boldsymbol{S}\\
& \boldsymbol{u}=\left(  \mathrm{u}_{n}\right)  _{n\in\mathbb{N}} &
\rightarrow & \boldsymbol{v}=\left(  \mathrm{v}_{n}\right)  _{n\in\mathbb{Z}}%
\end{array}
,
\]
where $\boldsymbol{v}=\left(  \mathrm{v}_{n}\right)  _{n\in\mathbb{Z}}$ is the
unique vector of $\boldsymbol{S}$ satisfying $\mathrm{v}_{n}=\mathrm{u}_{-n}$
for all $n\leq0$. The vector $\Theta(\boldsymbol{u})\in\boldsymbol{S}$ is
called the solution of the $Fib_{p}$ recurrence for the initial condition
$\boldsymbol{u}\mathbf{\in}\boldsymbol{U}$. The vector space $\boldsymbol{U}$
is called the space of initial conditions.

In order to analyze the asymptotic behavior of a solution
\[
\Theta(\boldsymbol{u})=\left(  \mathrm{v}_{n}\right)  _{n\in\mathbb{Z}}%
\in\boldsymbol{S},
\]
we define the generating function $G(\boldsymbol{u})$ as the formal power
series with coefficients in $\mathbb{C}^{p}$%
\[
G(\boldsymbol{u})=\sum\limits_{n\geq0}\mathrm{v}_{n}z^{n}.
\]
Alternatively, $G(\boldsymbol{u})$ can be defined as the element of the
$\mathbb{C}$-vector space $\mathbb{C}[[z]]^{p}$%
\[
G(\boldsymbol{u})=(G_{1}(\boldsymbol{u}),...,G_{p}(\boldsymbol{u}))\text{,}%
\]
with%
\begin{equation}
G_{\alpha}(\boldsymbol{u})=\sum\limits_{n\geq0}\mathrm{v}_{n}^{\left(
\alpha\right)  }z^{n}\in\mathbb{C}[[z]]\text{, }\alpha=1,\ldots,p\text{,}
\label{Galpha}%
\end{equation}
where $\mathrm{v}_{n}^{\left(  \alpha\right)  }$ denotes the $\alpha$-th
component of $\mathrm{v}_{n}$ with respect the standard base of $\mathbb{C}%
^{p}$.

Now we introduce the main ingredients of this work: the notions of
\textit{kneading matrix and kneading determinant} of a $Fib_{p}$ recurrence.
The idea is to look at the sequence $\left(  \mathbf{A}_{n}\right)
_{n\in\mathbb{N}}$ as a generating function $\sum_{n\geq0}\mathbf{A}_{n}z^{n}$
with coefficients in $\mathbb{C}^{p\times p}$. Naturally, this generating
function can be identified with the $p\times p$ matrix of formal power series%
\begin{equation}
\mathbf{K}=\left(
\begin{array}
[c]{ccc}%
K(1,1) & \cdots & K(1,p)\\
\vdots & \ddots & \vdots\\
K(p,1) & \cdots & K(p,p)
\end{array}
\right)  \text{,} \label{KneadMatrixK}%
\end{equation}
defined by%
\[
K(i,j)=%
{\displaystyle\sum\limits_{n\geq0}}
A_{n}(i,j)z^{n}\in\mathbb{C}[[z]]\text{.}%
\]
The matrix $\mathbf{K}$ is called the kneading matrix of the $Fib_{p}$
recurrence and the invertible formal power series%
\begin{equation}
\Delta=\det(\mathbf{I}-z\mathbf{K}) \label{Kneadingdet}%
\end{equation}
is called the kneading determinant of the $Fib_{p}$ recurrence.

Trivially, the entries of the kneading matrix are polynomials if and only if
the $Fib_{p}$ recurrence is of finite order. Hence, the kneading determinant
of a $Fib_{p}^{k}$ recurrence is actually a polynomial.

\section{Main results}

The first result of this work concerns the particular case of linear
recurrences of finite order and shows that the definition of kneading
determinant extends the usual definition of discriminant of a finite order
linear recurrence as defined in classical textbooks \cite{rabbit,saber}.

Recall that the Frobenius companion matrix\textbf{\ }of a $Fib_{p}^{k}$
recurrence is the $kp\times kp$ matrix%

\begin{equation}
\mathbf{F}=\left(
\begin{array}
[c]{cccc}%
\mathbf{A}_{0} & \cdots & \mathbf{A}_{k-2} & \mathbf{A}_{k-1}\\
\mathbf{I} & \cdots & \mathbf{0} & \mathbf{0}\\
\vdots & \ddots & \vdots & \vdots\\
\mathbf{0} & \cdots & \mathbf{I} & \mathbf{0}%
\end{array}
\right)  , \label{Compan}%
\end{equation}
where $\mathbf{I}$ and $\mathbf{0}$ denote respectively the $p\times p$
identity matrix and the $p\times p$ zero matrix. As the next result shows, the
classical discriminant $\det(\mathbf{I}-z\mathbf{F})$ coincides with the
kneading determinant of the $Fib_{p}^{k}$ recurrence.

\begin{theorem}
\label{T1}For any $Fib_{p}^{k}$ recurrence the relation%
\[
\det(\mathbf{I}-z\mathbf{F})=\det(\mathbf{I}-z\mathbf{K})
\]
holds.
\end{theorem}

The previous theorem is useful to compute explicitly the discriminant of a
vectorial finite recurrence $Fib_{p}^{k}$. A simple example illustrates this feature.

\begin{example}
\label{secondexample}Except for the case $p=1$, the computation of the
discriminant $\det(\mathbf{I}-z\mathbf{F})$ by standard methods requires in
general a large number of tedious computations. As an example, consider the
$Fib_{2}^{3}$ recurrence defined by%
\[
\mathbf{A}_{0}=\left(
\begin{array}
[c]{rr}%
1 & \quad1\\
1 & 1
\end{array}
\right)  \text{, }\mathbf{A}_{1}=\left(
\begin{array}
[c]{cc}%
-1 & -1\\
-1 & -1
\end{array}
\right)  ,\text{ }\mathbf{A}_{2}=\left(
\begin{array}
[c]{rr}%
0 & \quad1\\
1 & 0
\end{array}
\right)  \text{.}%
\]
The companion matrix is%
\[
\mathbf{F}=\left(
\begin{array}
[c]{rrrrrr}%
1 & \quad1 & -1 & -1 & \quad0 & \quad1\\
1 & 1 & -1 & -1 & 1 & 0\\
1 & 0 & \quad0 & \quad0 & 0 & 0\\
0 & 1 & 0 & 0 & 0 & 0\\
0 & 0 & 1 & 0 & 0 & 0\\
0 & 0 & 0 & 1 & 0 & 0
\end{array}
\right)  .
\]
After some cumbersome computations, one gets%
\[
\det(\mathbf{I}-z\mathbf{F})=\left(  1-z\right)  \left(  1+z\right)  \left(
1-z+z^{2}\right)  ^{2}.
\]
On the other hand, as the kneading matrix is%
\[
\mathbf{K}=\left(
\begin{array}
[c]{cc}%
1-z & 1-z+z^{2}\\
1-z+z^{2} & 1-z
\end{array}
\right)  ,
\]
a simple computation gives%
\[
\Delta=\det(\mathbf{I}-z\mathbf{K})=\left(  1-z\right)  \left(  1+z\right)
\left(  1-z+z^{2}\right)  ^{2}.
\]
which agrees with the value of $\det(\mathbf{I}-z\mathbf{F})$ obtained above
by direct approach.
\end{example}

Next, we focus on the main topic of this work: the study of the asymptotic
behavior of the solutions of a $Fib_{p}$ recurrence. Our first goal is to
provide explicit formulas for computing the generating functions
$G(\boldsymbol{u})$ of a finite or infinite order $Fib_{p}$ recurrence.

First of all observe that from the linearity of the map
\[%
\begin{array}
[c]{cccc}%
G & :\boldsymbol{U} & \rightarrow & \mathbb{C}[[z]]^{p}\\
& \boldsymbol{u} & \rightarrow & G(\boldsymbol{u}),
\end{array}
\]
one has%

\begin{equation}
G(\boldsymbol{u})=\sum\limits_{\beta\geq1}c_{\beta}G(\boldsymbol{e}_{\beta
})=\sum\limits_{\beta\geq1}c_{\beta}\left(  G_{1}(\boldsymbol{e}_{\beta
}),...,G_{p}(\boldsymbol{e}_{\beta})\right)  \text{\textbf{,}} \label{FP12}%
\end{equation}
where $\left(  c_{\beta}\right)  _{\beta\in\mathbb{Z}^{+}}$ denotes the
coordinates of $\boldsymbol{u}$ with respect to the standard basis $\left(
\boldsymbol{e}_{\beta}\right)  _{\beta\in\mathbb{Z}^{+}}$ of $\boldsymbol{U}$.
Therefore, to accomplish this task we just need to focus on the generating
functions $G_{\alpha}(\boldsymbol{e}_{\beta})$.

For this purpose, we define for each $\alpha=1,...,p$ and each $\beta
\in\mathbb{Z}^{+}$ the extended kneading matrix $\mathbf{K}_{\alpha}(\beta)$
adding one more row and one more column to the kneading matrix $\mathbf{K}$ of
the $Fib_{p}$ recurrence. More precisely we define $\mathbf{K}_{\alpha}%
(\beta)\in\mathbb{C}[[z]]^{\left(  p+1\right)  \times\left(  p+1\right)  }$ by
setting
\begin{equation}
\mathbf{K}_{\alpha}(\beta)=\left(
\begin{array}
[c]{cccc}%
K\left(  1,1\right)  & \cdots & K\left(  1,p\right)  & K\left(  1,\beta\right)
\\
\vdots & \ddots & \vdots & \vdots\\
K\left(  p,1\right)  & \cdots & K\left(  p,p\right)  & K\left(  p,\beta\right)
\\
\delta(\alpha,1) & \cdots & \delta(\alpha,p) & \delta(\alpha,\beta)
\end{array}
\right)  , \label{KneadMatrixKalfa}%
\end{equation}
where $\delta(i,j)$ is the usual Kronecker delta function. For the last column
of $\mathbf{K}_{\alpha}(\beta)$\ we consider the quotient $q$ and the reminder
$r$ of the division of $\beta$ by $p$ to introduce%
\[
K\left(  i,\beta\right)  =\left\{
\begin{array}
[c]{ll}%
\underset{n\geq0}{\sum}A_{n+q-1}(i,p)z^{n},\text{ } & \text{if }p\text{
divides }\beta\\
\underset{n\geq0}{\sum}A_{n+q}(i,r)z^{n},\text{ } & \text{otherwise.}%
\end{array}
\right.
\]
Finally, we define the extended kneading determinant%
\[
\Delta_{\alpha}(\beta)=\det(\mathbf{I}-z\mathbf{K}_{\alpha}(\beta))\text{.}%
\]

Now, we can state the main and new result of this work which gives explicitly
the entries of matrix generating function $G\left(  z\right)  $ for the
solutions of any vectorial recurrence. To our knowledge there is no other way
of computing explicitly the solutions of any $Fib_{p}^{\infty}$ recurrence.
Naturally, this result solves also the classical problem of computing the
solutions of finite order linear recurrences, which is classically done using
Jordan canonical forms \cite{rabbit,saber}.

\begin{theorem}
\label{T2}For every $\alpha=1,..,p$ and every vector $\boldsymbol{e}_{\beta}$
of the standard basis of $\boldsymbol{U}$, the generating function $G_{\alpha
}(\boldsymbol{e}_{\beta})$ of a $Fib_{p}$ recurrence satisfies the following
equality in $\mathbb{C}[[z]]$%
\[
zG_{\alpha}(\boldsymbol{e}_{\beta})=1-\Delta^{-1}\Delta_{\alpha}%
(\beta)\text{.}%
\]

\end{theorem}

\begin{example}
\label{E1}In order to illustrate Theorem \ref{T2}, we compute the generating
functions $G(\boldsymbol{e}_{1})$ of the $Fib_{2}^{\infty}$ recurrence defined
by%
\[
\mathbf{A}_{n}=\frac{1}{n!}\left(
\begin{array}
[c]{cc}%
-\frac{1}{n+1}\quad\medskip & 2^{n}\medskip\\
\left(  -1\right)  ^{n}\quad & 0
\end{array}
\right)  \text{, for }n\in\mathbb{N}\text{.}%
\]
We have%
\[
\mathbf{K=}\left(
\begin{array}
[c]{cc}%
\frac{1-e^{z}}{z}\medskip & \quad e^{2z}\\
e^{-z} & 0
\end{array}
\right)  \text{ and }\Delta=\det\left(  \mathbf{I}-z\mathbf{K}\right)
=\left(  1-z^{2}\right)  e^{z}.
\]
On the other hand, as the extended kneading matrices $\mathbf{K}_{1}(1)$ and
$\mathbf{K}_{2}(1)$ are defined by:
\[
\mathbf{K}_{\alpha}(1)=\left(
\begin{array}
[c]{ccc}%
\frac{1-e^{z}}{z} & e^{2z} & \frac{1-e^{z}}{z}\\
e^{-z} & 0 & e^{-z}\\
\delta(\alpha,1) & \delta(\alpha,2) & \delta(\alpha,1)
\end{array}
\right)  \text{, for }\alpha=1,2\text{,}%
\]
the extended kneading determinants are%
\[
\mathbf{\Delta}_{1}(1)=\det\left(  \mathbf{I}-z\mathbf{K}_{1}(1)\right)
=\left(  1-z^{2}\right)  e^{z}-z
\]
and%
\[
\mathbf{\Delta}_{2}(1)=\det\left(  \mathbf{I}-z\mathbf{K}_{2}(1)\right)
=\left(  1-z^{2}\right)  e^{z}-z^{2}e^{-z}.
\]
Finally, by Theorem \ref{T2} we have
\[
G(\boldsymbol{e}_{1})=\left(  \frac{1}{\left(  1-z^{2}\right)  e^{z}},\frac
{z}{\left(  1-z^{2}\right)  e^{2z}}\right)  \text{.}%
\]

\end{example}

We complete this section by discussing the existence of asymptotic
closed-forms for the solutions of an infinite order linear recurrence. As we
will see, Theorem \ref{T2} plays a central role in this discussion.

As a motivation for Theorem \ref{T3}, the last result of this section, we
recall the case of finite order recurrences where the existence of closed
forms for the solutions are well known.

Let $\lambda_{1},...,\lambda_{m}\in\mathbb{C}\backslash\left\{  0\right\}  $
and $\mathrm{mul}\left(  \lambda_{1}\right)  ,...,\mathrm{mul}\left(
\lambda_{m}\right)  \in\mathbb{Z}^{+}$ denote the nonzero eigenvalues and
corresponding algebraic multiplicities of the companion matrix, $\mathbf{F}$,
of a finite recurrence $Fib_{p}^{k}$, then for any solution $\Theta
(\boldsymbol{u})=\left(  \mathrm{v}_{n}^{\left(  1\right)  },...,\mathrm{v}%
_{n}^{\left(  p\right)  }\right)  $ and every $\alpha=1,...,p$, there exist
unique constants $c_{i,j}^{\left(  \alpha\right)  }(\boldsymbol{u}%
)\in\mathbb{C}$, with $i=1,...,m$ and $j=1,...,\mathrm{mul}\left(  \lambda
_{i}\right)  $, such that%
\[
\mathrm{v}_{n}^{\left(  \alpha\right)  }=\sum_{i=1}^{m}\sum_{j=1}%
^{\mathrm{mul}(\lambda_{i})}\frac{c_{i,j}^{\left(  \alpha\right)
}(\boldsymbol{u})(n+j-1)!}{\left(  j-1\right)  !n!}\lambda_{i}^{n}\text{, for
all }n>kp.
\]
In particular, if $\mathrm{mul}\left(  \lambda_{i}\right)  =1$, for
$i=1,...,p$, one gets the Binet formula%
\[
\mathrm{v}_{n}^{\left(  \alpha\right)  }=\sum_{i=1}^{m}c_{i,i}^{\left(
\alpha\right)  }\lambda_{i}^{n}\text{, for all }n>kp.
\]
For the original Binet Formula and historical approach see page 281 of
\cite{Burton}.

In the case of infinite order recurrences with nonrational\footnote{A formal
power series $a=\sum_{n\geq0}a_{n}z^{n}$ is said to be a rational function of
$z,$ if there exist polynomials $p,q$ such that $q$ is invertible in
$\mathbb{C}\left[  [z]\right]  $ and $a=q^{-1}p$.} generating functions as
seen in Example \ref{E1}, there are no closed-forms. However, one can
establish the existence of asymptotic closed forms in some cases. For
instance, in \cite{rash1} are obtained Binet formulas for periodic
$F_{1}^{\infty}$ recurrences.

In the case of infinite vector recurrences of the type $F_{p}^{\infty}$ we
need to introduce some essential concepts to state Theorem \ref{T3}.

As usual, a matrix $\mathbf{M\in}\mathbb{C}[[z]]^{p\times p}$ is said to be
holomorphic\footnote{A formal power series $a=\sum_{n\geq0}a_{n}z^{n}$ with
radius of convergence $R$ is said to be holomorphic on $D_{\rho}$ if
$R\geq\rho$. Similarly, one says that $a\in\mathbb{C}\left[  [z]\right]  $ is
meromorphic on $D_{\rho}$ if there exist $b,c\in\mathbb{C}\left[  [z]\right]
$ such that $b$ and $c$ are holomorphic on $D_{\rho}$, $c$ is invertible in
$\mathbb{C}\left[  [z]\right]  $ and $a=c^{-1}b$.} (resp. meromorphic) on the
open disk $D_{\rho}=\left\{  z\in\mathbb{C}:\left\vert z\right\vert
<\rho\right\}  $, with $\rho\in\left]  0,+\infty\right]  $, if the entries of
$\mathbf{M}$ are holomorphic (resp. meromorphic) functions on $D_{\rho}$.

Consequently, if the kneading matrix $\mathbf{K}$ of a $Fib_{p}$ recurrence is
holomorphic on $D_{\rho}$, we can look at the kneading determinant $\Delta$ as
an analytic function on $D_{\rho}$.

The next definition is motivated by Theorem \ref{T1}, which proves that
$\lambda\in\mathbb{C}\backslash\left\{  0\right\}  $ is an eigenvalue of the
companion matrix of a finite recurrence $Fib_{p}^{k}$ if and only if the
kneading determinant $\Delta$ has a zero at $\lambda^{-1}$.

\begin{definition}
Assume that the kneading matrix, $\mathbf{K}$, of a $Fib_{p}$-recurrence is
holomorphic on $D_{\rho}$. A complex number, $\lambda$, with $\left\vert
\lambda\right\vert >\rho^{-1}$, is said to be a generalized eigenvalue of the
$Fib_{p}$-recurrence with multiplicity $\mathrm{mul}\left(  \lambda\right)
\in\mathbb{Z}^{+}$ if the kneading determinant $\Delta$ has a $\mathrm{mul}%
\left(  \lambda\right)  $-order zero at $\lambda^{-1}$. A generalized
eigenvalue, $\lambda$, is said to be dominant if $\left\vert \lambda
\right\vert \geq1$.
\end{definition}

Notice that if the kneading matrix $\mathbf{K}$ is holomorphic on some
$D_{\rho}$, with $\rho>1$, then the $Fib_{p}$-recurrence has finitely many
generalized eigenvalues. This is the setting of the second main result of the paper.

\begin{theorem}
\label{T3}Let $\lambda_{1},...,\lambda_{m}$ be the dominant eigenvalues of a
$Fib_{p}$ recurrence whose kneading matrix $K$ is holomorphic on some open
disk $D_{\rho}$, with $\rho>1$. Then, for any solution $\Theta(\boldsymbol{u}%
)=\left(  \mathrm{v}_{n}^{\left(  1\right)  },...,\mathrm{v}_{n}^{\left(
p\right)  }\right)  $ and every $\alpha=1,...,p$ there exist unique constants
$c_{i,j}^{\left(  \alpha\right)  }(\boldsymbol{u})\in\mathbb{C}$, with
$i=1,...,m$ and $j=1,...,\mathrm{mul}\left(  \lambda_{i}\right)  $ such that%
\[
\underset{n\rightarrow+\infty}{\lim}\left(  \mathrm{v}_{n}^{\left(
\alpha\right)  }-\sum_{i=1}^{m}\sum_{j=1}^{\mathrm{mul}(\lambda_{i})}%
\frac{c_{i,j}^{\left(  \alpha\right)  }(\boldsymbol{u})(n+j-1)!}{\left(
j-1\right)  !n!}\lambda_{i}^{n}\right)  =0.
\]

\end{theorem}

With this theorem it is clear that the dominant eigenvalues of an infinite
linear recurrence characterize the asymptotic behavior of the solutions of
that recurrence.

Note that the previous theorem \ shows that if $\mathrm{mul}(\lambda_{i})=1,$
$i=1,...,m,$ then%
\[
\underset{n\rightarrow+\infty}{\lim}\left(  \mathrm{v}_{n}^{\left(
\alpha\right)  }-\sum_{i=1}^{m}d_{i}^{\left(  \alpha\right)  }(\boldsymbol{u}%
)\lambda_{i}^{n}\right)  =0\text{, with }d_{i}^{\left(  \alpha\right)
}=c_{i,i}^{\left(  \alpha\right)  }\text{,}%
\]
which is the generalization of the classic Binet formula.

We finish this section with an example illustrating Theorem \ref{T3}.

\begin{example}
Let us return to Example \ref{E1}. From Theorem \ref{T3} it is easy to prove
that any solution $\Theta(\boldsymbol{u})=\left(  \mathrm{v}_{n}^{\left(
1\right)  },\mathrm{v}_{n}^{\left(  2\right)  }\right)  $ is asymptotically
periodic with period $2$, that is both sequences $\left(  \mathrm{v}%
_{2n}^{\left(  1\right)  },\mathrm{v}_{2n}^{\left(  2\right)  }\right)  $ and
$\left(  \mathrm{v}_{2n+1}^{\left(  1\right)  },\mathrm{v}_{2n+1}^{\left(
2\right)  }\right)  $ are convergent. Indeed, since the kneading matrix is
holomorphic on $\mathbb{C}$ and $\Delta=(1-z^{2})e^{z}$, the dominant
eigenvalues are $\lambda_{1}=1$ and $\lambda_{2}=-1$, with $\mathrm{mul}%
\left(  \lambda_{1}\right)  =\mathrm{mul}\left(  \lambda_{2}\right)  =1$.
Therefore, for any solution $\Theta(\boldsymbol{u})$ there exist unique
constants $c_{1}^{\left(  1\right)  }(\boldsymbol{u}),c_{2}^{\left(  1\right)
}(\boldsymbol{u}),c_{1}^{\left(  2\right)  }(\boldsymbol{u}),c_{2}^{\left(
2\right)  }(\boldsymbol{u})\in\mathbb{C}$ such that%
\[
\underset{n\rightarrow+\infty}{\lim}\left(  \mathrm{v}_{n}^{\left(  1\right)
}-c_{1}^{\left(  1\right)  }(\boldsymbol{u})-c_{2}^{\left(  1\right)
}(\boldsymbol{u})(-1)^{n}\right)  =0
\]
and%
\[
\underset{n\rightarrow+\infty}{\lim}\left(  \mathrm{v}_{n}^{\left(  2\right)
}-c_{1}^{\left(  2\right)  }(\boldsymbol{u})-c_{2}^{\left(  2\right)
}(\boldsymbol{u})(-1)^{n}\right)  =0\text{.}%
\]
Hence, we have%
\[
\underset{n\rightarrow+\infty}{\lim}\left(  \mathrm{v}_{2n}^{\left(  1\right)
},\mathrm{v}_{2n}^{\left(  2\right)  }\right)  =\left(  c_{1}^{\left(
1\right)  }(\boldsymbol{u})+c_{2}^{\left(  1\right)  }(\boldsymbol{u}%
),c_{1}^{\left(  2\right)  }(\boldsymbol{u})+c_{2}^{\left(  2\right)
}(\boldsymbol{u})\right)
\]
and%
\[
\underset{n\rightarrow+\infty}{\lim}\left(  \mathrm{v}_{2n+1}^{\left(
1\right)  },\mathrm{v}_{2n+1}^{\left(  2\right)  }\right)  =\left(
c_{1}^{\left(  1\right)  }(\boldsymbol{u})-c_{2}^{\left(  1\right)
}(\boldsymbol{u}),c_{1}^{\left(  2\right)  }(\boldsymbol{u})-c_{2}^{\left(
2\right)  }(\boldsymbol{u})\right)  .
\]

\end{example}

The rest of the paper will be devoted to the proofs of Theorems \ref{T1},
\ref{T2} and \ref{T3}.

\section{Pairs of linear endomorphisms\label{S3}}

The proofs of Theorems \ref{T1} and \ref{T2} on the next section are rooted in
the main Theorem of \cite{jalves1} concerning the determinant of a pair of
linear endomorphisms with finite rank. This last theorem extends to a wider
context the well known relationship between discriminant and traces for a
matrix $\mathbf{X}\in\mathbb{C}^{m\times m}$%
\begin{equation}
\det(\mathbf{I}-z\boldsymbol{X})=\exp\sum_{n\geq1}-\frac{tr(\mathbf{X}^{n}%
)}{n}z^{n}. \label{FP1}%
\end{equation}

In order to improve the readability of the paper we present a brief
description of the results obtained in \cite{jalves1}.

Throughout this section, $U$ denotes an arbitrary (finite or infinite
dimensional) vector space over $\mathbb{C}$; the space of linear forms on $U$
will be denoted by $U^{\ast}$ and the space of linear endomorphisms on $U$
will be denoted by $\mathrm{L}(U)$. If $\psi\in\mathrm{L}(U)$ and $n$ is a
nonnegative integer, the $n$-th iterate $\psi^{n}$ is defined recursively by
$\psi^{0}=Id_{U}\in\mathrm{L}(U)$, $\psi^{n}=\psi\circ\psi^{n-1}\in
\mathrm{L}(U)$, for $n\geq1$.

Recall that a linear endomorphism $\psi\in\mathrm{L}(U)$ is said to have
finite rank if there exist vectors $u_{1},...,u_{p}\in U$ and linear forms
$\omega_{1},...,\omega_{p}\in U^{\ast}$ such that
\[
\psi=\omega_{1}\otimes u_{1}+\omega_{2}\otimes u_{2}+\cdots+\omega_{p}\otimes
u_{p},
\]
with the usual notation
\[
\omega\in U^{\ast},\text{ }u\in U:\left(  \omega\otimes u\right)  \left(
x\right)  =\omega\left(  x\right)  u,\;x\in U.
\]
The subspace of $\mathrm{L}(U)$ whose elements are the linear endomorphisms
with finite rank on $U$ will be denoted by $\mathrm{L}_{{\small FR}}(U)$.

The importance of $\mathrm{L}_{{\small FR}}(U)$, in this context, lies in the
existence of the trace for any $\psi\in\mathrm{L}_{{\small FR}}(U)$, trace
that is not evidently defined for an arbitrary $\psi\in\mathrm{L}(U)$.

Let us then introduce the following definition.

\begin{definition}
A pair of endomorphisms $\left(  \varphi,\psi\right)  \in\mathrm{L}%
(U)\times\mathrm{L}(U)$ is said to have finite rank if $\psi-\varphi
\in\mathrm{L}_{{\small FR}}(U)$.
\end{definition}

Notice that if a pair $\left(  \varphi,\psi\right)  $ has finite rank, then
the pair $\left(  \varphi^{n},\psi^{n}\right)  $ has finite rank for all
$n\geq0$. Therefore, the trace of $\varphi^{n}-\psi^{n}$ is defined and the
following definition makes sense.

\begin{definition}
For any pair $\left(  \varphi,\psi\right)  \in\mathrm{L}(U)\times
\mathrm{L}(U)$ with finite rank, the determinant of $\left(  \varphi
,\psi\right)  $ is defined as the formal power series%
\[
\Delta\left(  \varphi,\psi\right)  =\exp\sum_{n\geq1}\frac{\mathrm{tr}%
(\varphi^{n}-\psi^{n})}{n}z^{n}\text{.}%
\]

\end{definition}

If a pair $\left(  \varphi,\psi\right)  $ has finite rank, then $\left(
\psi,\varphi\right)  $ has finite rank too and
\[
\Delta\left(  \varphi,\psi\right)  \Delta\left(  \psi,\varphi\right)  =1.
\]
Thus, $\Delta\left(  \varphi,\psi\right)  $ is invertible in $\mathbb{C}[[z]]$
with inverse%
\[
\left[  \Delta\left(  \varphi,\psi\right)  \right]  ^{-1}=\Delta\left(
\psi,\varphi\right)  \text{.}%
\]
More generally one has the following proposition.

\begin{proposition}
\label{PP1}If $\left(  \varphi,\psi\right)  \in\mathrm{L}(U)\times
\mathrm{L}(U)$ and $\left(  \psi,\chi\right)  \in\mathrm{L}(U)\times
\mathrm{L}(U)$ have both finite rank, then $\left(  \varphi,\chi\right)
\in\mathrm{L}(U)\times\mathrm{L}(U)$ has finite rank and $\Delta\left(
\varphi,\chi\right)  =$ $\Delta\left(  \varphi,\psi\right)  \Delta\left(
\psi,\chi\right)  $.
\end{proposition}

Notice that if the space $U$ is finite dimensional, then every pair $\left(
\varphi,\psi\right)  $ has finite rank and by (\ref{FP1}) one gets%
\[
\Delta\left(  \varphi,\psi\right)  =\frac{\det(\mathbf{I}-z\mathbf{Y})}%
{\det(\mathbf{I}-z\mathbf{X})},
\]
where $\mathbf{X}$ (respectively $\mathbf{Y}$) is the matrix that represents
$\varphi$ (respectively $\psi$) with respect to some basis of $U$. So, in this
particular case $\Delta\left(  \varphi,\psi\right)  $ is a rational function
of $z$.

The situation becomes entirely different if the space $U$ is infinite
dimensional. In this case the rationality of $\Delta\left(  \varphi
,\psi\right)  $ fails in general. This fact is a simple consequence of the
next result which enables us to express $\Delta\left(  \varphi,\psi\right)  $
in terms of determinants. To state it we have to introduce some additional notation.

Observe that if a pair $\left(  \varphi,\psi\right)  \in\mathrm{L}(U)\times$
$\mathrm{L}(U)$ has finite rank, then there exist vectors $u_{1},...,u_{p}\in
U$ and linear forms $\omega_{1},...,\omega_{p}\in U^{\ast}$ such that
\begin{equation}
\psi-\varphi=\omega_{1}\otimes u_{1}+\omega_{2}\otimes u_{2}+\cdots+\omega
_{p}\otimes u_{p} \label{FP2}%
\end{equation}
and so, we can define the matrix $\mathbf{M}\in\mathbb{C}\left[  \left[
z\right]  \right]  ^{p\times p}$ by setting
\begin{equation}
\mathbf{M}=\left(
\begin{array}
[c]{ccc}%
\underset{n\geq0}{\sum}\omega_{1}\varphi^{n}(u_{1})z^{n} & \cdots &
\underset{n\geq0}{\sum}\omega_{1}\varphi^{n}(u_{p})z^{n}\\
\vdots & \ddots & \vdots\\
\underset{n\geq0}{\sum}\omega_{p}\varphi^{n}(u_{1})z^{n} & \cdots &
\underset{n\geq0}{\sum}\omega_{p}\varphi^{n}(u_{p})z^{n}%
\end{array}
\right)  \text{.} \label{FP3}%
\end{equation}
Now we can state the main Theorem of \cite{jalves1} which establishes a
fundamental relationship between $\Delta\left(  \varphi,\psi\right)  $ and the
determinant of the matrix $\mathbf{I}-z\mathbf{M}$.

\begin{lemma}
\label{LP1}Let $\left(  \varphi,\psi\right)  \in\mathrm{L}(U)\times
\mathrm{L}(U)$ be a pair of endomorphisms with finite rank. If the vectors
$u_{1},...,u_{p}\in U$ and the linear forms $\omega_{1},...,\omega_{p}\in
U^{\ast}$ satisfy (\ref{FP2}), then $\Delta\left(  \varphi,\psi\right)
=\det(\mathbf{I}-z\mathbf{M})$.
\end{lemma}

Two consequences of this result are needed.

The first one can be thought as an alternative method for computing the
discriminant $\det(\mathbf{I}-z\mathbf{Y})$ of a complex matrix $\mathbf{Y}%
\in\mathbb{C}^{m\times m}$.

The idea is to consider a nilpotent matrix $\mathbf{X}\in\mathbb{C}^{m\times
m}$ and to look at $(\mathbf{X},\mathbf{Y})$ as a pair of linear endomorphism
on $\mathbb{C}^{m}$ with finite rank. So, we can consider column matrices
$\mathbf{C}_{1},...,\mathbf{C}_{p}\in\mathbb{C}^{m\times1}$ and row matrices
$\mathbf{R}_{1},...,\mathbf{R}_{p}\in\mathbb{C}^{1\times m}$ satisfying%
\begin{equation}
\mathbf{Y}-\mathbf{X}=\mathbf{C}_{1}\mathbf{R}_{1}+\mathbf{C}_{2}%
\mathbf{R}_{2}+\cdots+\mathbf{C}_{p}\mathbf{R}_{p} \label{FP4}%
\end{equation}
and by Lemma \ref{LP1} we can write%
\begin{equation}
\exp\sum_{n\geq1}\frac{\mathrm{tr}(\mathbf{X}^{n}-\mathbf{Y}^{n})}{n}%
z^{n}=\det(\mathbf{I}-z\mathbf{M})\text{,} \label{FP8}%
\end{equation}
with%
\begin{equation}
\mathbf{M}=\left(
\begin{array}
[c]{ccc}%
\underset{n\geq0}{\sum}\mathbf{R}_{1}\mathbf{X}^{n}\mathbf{C}_{1}z^{n} &
\cdots & \underset{n\geq0}{\sum}\mathbf{R}_{1}\mathbf{X}^{n}\mathbf{C}%
_{p}z^{n}\\
\vdots & \ddots & \vdots\\
\underset{n\geq0}{\sum}\mathbf{R}_{p}\mathbf{X}^{n}\mathbf{C}_{1}z^{n} &
\cdots & \underset{n\geq0}{\sum}\mathbf{R}_{p}\mathbf{X}^{n}\mathbf{C}%
_{p}z^{n}%
\end{array}
\right)  . \label{FP5}%
\end{equation}
Moreover, being $\mathbf{X}$ nilpotent then $\mathbf{M}$ is a $p\times p$
matrix of polynomials. Hence $\det(\mathbf{I}-z\mathbf{M})$ is a polynomial
too. As we will see in the next result, this polynomial is actually the
discriminant of $\mathbf{Y}$.

\begin{theorem}
\label{LP2}Let $\mathbf{X}\in\mathbb{C}^{m\times m},\mathbf{Y}\in
\mathbb{C}^{m\times m},\mathbf{C}_{1},...,\mathbf{C}_{p}\in\mathbb{C}%
^{m\times1}$ and $\mathbf{R}_{1},...,\mathbf{R}_{p}\in\mathbb{C}^{1\times m}$
satisfying (\ref{FP4}). If $\mathbf{X}$ is nilpotent, then the equality
$\det(\mathbf{I}-z\mathbf{Y})=\det(\mathbf{I}-z\mathbf{M})$ holds in
$\mathbb{C[}z]$.
\end{theorem}

\begin{pf}
As $\mathbf{X}$ is nilpotent one has $tr(\mathbf{X}^{n})=0$ for $n\geq1$.
Combining this with (\ref{FP1}) and (\ref{FP8}) one gets%
\begin{align*}
\det(\mathbf{I}-z\mathbf{Y})  &  =\exp\sum_{n\geq1}-\frac{\mathrm{tr}%
(\mathbf{Y}^{n})}{n}z^{n}\\
&  =\exp\sum_{n\geq1}\frac{\mathrm{tr}(\mathbf{X}^{n}-\mathbf{Y}^{n})}{n}%
z^{n}\\
&  =\det(\mathbf{I}-z\mathbf{M}),
\end{align*}
as desired.
\end{pf}

A second consequence of Lemma \ref{LP1} concerns the general and difficult
problem of studying the analytic properties of the generating function%
\begin{equation}
\underset{n\geq0}{\sum}\omega\psi^{n}(u)z^{n}\text{,} \label{FP6}%
\end{equation}
where $\psi\in\mathrm{L}(U)$, $\omega\in U^{\ast}$ and $u\in U$ are arbitrary.

An idea that can be useful, is to consider a pair $\left(  \varphi
,\psi\right)  \in\mathrm{L}(U)\times\mathrm{L}(U)$ with finite rank and write
(\ref{FP6}) in terms of determinants with the desired analytic properties.

Notice that if a pair $\left(  \varphi,\psi\right)  $ has finite rank, then
$\left(  \varphi,\psi+\omega\otimes u\right)  $ has finite rank too. In fact,
if the vectors $u_{1},...,u_{p}\in U$ and the linear forms $\omega
_{1},...,\omega_{p}\in U^{\ast}$ satisfy (\ref{FP2}), then%
\[
\left(  \psi+\omega\otimes u\right)  -\varphi=\omega_{1}\otimes u_{1}%
+\omega_{2}\otimes u_{2}+\cdots+\omega_{p}\otimes u_{p}+\omega\otimes u
\]
and by Lemma \ref{LP1}%
\begin{equation}
\Delta\left(  \varphi,\psi+\omega\otimes u\right)  =\det(\mathbf{I}%
-z\mathbf{M}_{\omega}(u)), \label{FP9}%
\end{equation}
where $\mathbf{M}_{\omega}(u)\in\mathbb{C}[[z]]^{(p+1)\times(p+1)}$ is the
extended matrix defined by%
\begin{equation}
\mathbf{M}_{\omega}(u)=\left(
\begin{array}
[c]{cccc}%
\underset{n\geq0}{\sum}\omega_{1}\varphi^{n}(u_{1})z^{n} & \cdots &
\underset{n\geq0}{\sum}\omega_{1}\varphi^{n}(u_{p})z^{n} & \underset{n\geq
0}{\sum}\omega_{1}\varphi^{n}(u)z^{n}\\
\vdots & \ddots & \vdots & \vdots\\
\underset{n\geq0}{\sum}\omega_{p}\varphi^{n}(u_{1})z^{n} & \cdots &
\underset{n\geq0}{\sum}\omega_{p}\varphi^{n}(u_{p})z^{n} & \underset{n\geq
0}{\sum}\omega_{p}\varphi^{n}(u)z^{n}\\
\underset{n\geq0}{\sum}\omega\varphi^{n}(u_{1})z^{n} & \cdots & \underset
{n\geq0}{\sum}\omega\varphi^{n}(u_{p})z^{n} & \underset{n\geq0}{\sum}%
\omega\varphi^{n}(u)z^{n}%
\end{array}
\right)  . \label{FP7}%
\end{equation}
Now, it is easy to establish a simple relationship between the generating
function of (\ref{FP6}) and the matrices $\mathbf{M}$ and $\mathbf{M}_{\omega
}(u)$ of (\ref{FP3}) and (\ref{FP7}).

\begin{lemma}
\label{LP3}Let $\left(  \varphi,\psi\right)  \in\mathrm{L}(U)\times
\mathrm{L}(U)$ be a pair with finite rank, $u\in U$ and $\omega\in U^{\ast}$.
If the vectors $u_{1},...,u_{p}\in U$ and the linear forms $\omega
_{1},...,\omega_{p}\in U^{\ast}$ satisfy (\ref{FP2}), then we have the
equality%
\[
\underset{n\geq0}{z\sum}\omega\psi^{n}(u)z^{n}=1-\frac{\det(\mathbf{I}%
-z\mathbf{M}_{\omega}(u))}{\det(\mathbf{I}-z\mathbf{M})}.
\]

\end{lemma}

\begin{pf}
Combining Lemma \ref{LP1} with (\ref{FP9}) and Proposition \ref{PP1}, one gets%
\begin{align*}
\frac{\det(\mathbf{I}-z\mathbf{M}_{\omega}(u))}{\det(\mathbf{I}-z\mathbf{M})}
&  =\frac{\Delta\left(  \varphi,\psi+\omega\otimes u\right)  }{\Delta\left(
\varphi,\psi\right)  }\\
&  =\Delta\left(  \psi,\varphi\right)  \Delta\left(  \varphi,\psi
+\omega\otimes u\right) \\
&  =\Delta\left(  \psi,\psi+\omega\otimes u\right)  .
\end{align*}
But, again by Lemma \ref{LP1} and because $\left(  \psi+\omega\otimes
u\right)  -\psi=\omega\otimes u$ we can write%
\[
\Delta\left(  \psi,\psi+\omega\otimes u\right)  =1-\underset{n\geq0}{z\sum
}\omega\psi^{n}(u)z^{n}\text{.}%
\]
Hence%
\[
1-\underset{n\geq0}{z\sum}\omega\psi^{n}(u)z^{n}=\frac{\det(\mathbf{I}%
-z\mathbf{M}_{\omega}(u))}{\det(\mathbf{I}-z\mathbf{M})}\text{,}%
\]
as desired.
\end{pf}

\section{Proofs of the main results}

At this stage we have all the ingredients to prove the main results of this
article: theorems \ref{T1}, \ref{T2} and \ref{T3}.

Theorem \ref{T1} is a simple consequence of Theorem \ref{LP2} given in the
previous section.

\begin{pf}
[Proof of Theorem $1$]Let $\mathbf{F\in}\mathbb{C}^{kp\times kp}$ be the
Frobenius companion matrix of a $Fib_{p}^{k}$ recurrence as defined in
(\ref{Compan}). For each $i=1,...,p$, let $\mathbf{R}_{i}\in\mathbb{C}%
^{1\times kp}$ be the $i$-th row of $\mathbf{F}$ and $\mathbf{C}_{i}%
\in\mathbb{C}^{kp\times1}$ the $i$-th vector of the standard basis of
$\mathbb{C}^{kp\times1}$. Evidently, the $kp\times kp$ matrix
\[
\mathbf{X}=\left(
\begin{array}
[c]{cccc}%
\mathbf{0} & \cdots & \mathbf{0} & \mathbf{0}\\
\mathbf{I} & \cdots & \mathbf{0} & \mathbf{0}\\
\vdots & \ddots & \vdots & \vdots\\
\mathbf{0} & \cdots & \mathbf{I} & \mathbf{0}%
\end{array}
\right)
\]
is nilpotent and%
\[
\mathbf{F}-\mathbf{X}=\mathbf{C}_{1}\mathbf{R}_{1}+\cdots+\mathbf{C}%
_{p}\mathbf{R}_{p}\text{.}%
\]
So, as the matrices $\mathbf{X},\mathbf{F},\mathbf{C}_{1},...,\mathbf{C}%
_{p},\mathbf{R}_{1},...,\mathbf{R}_{p}$ satisfy the assumptions of Theorem
\ref{LP2} we can write%
\begin{equation}
\det(\mathbf{I}-z\mathbf{F})=\det(\mathbf{I}-z\mathbf{M}), \label{FP15}%
\end{equation}
with%
\[
\mathbf{M}=\left(
\begin{array}
[c]{ccc}%
\underset{n\geq0}{\sum}\mathbf{R}_{1}\mathbf{X}^{n}\mathbf{C}_{1}z^{n} &
\cdots & \underset{n\geq0}{\sum}\mathbf{R}_{1}\mathbf{X}^{n}\mathbf{C}%
_{p}z^{n}\\
\vdots & \ddots & \vdots\\
\underset{n\geq0}{\sum}\mathbf{R}_{p}\mathbf{X}^{n}\mathbf{C}_{1}z^{n} &
\cdots & \underset{n\geq0}{\sum}\mathbf{R}_{p}\mathbf{X}^{n}\mathbf{C}%
_{p}z^{n}%
\end{array}
\right)  .
\]
But by (\ref{KneadMatrixK}) and because $\mathbf{A}_{n}=\mathbf{0}$ for $n\geq
k$, the $(i,j)$ entry of $\mathbf{M}$ is%
\[
\underset{n\geq0}{\sum}\mathbf{R}_{i}\mathbf{X}^{n}\mathbf{C}_{j}%
z^{n}=\underset{n=0}{\sum^{k-1}}A_{n}(i,j)z^{n}=\underset{n\geq0}{\sum}%
A_{n}(i,j)z^{n}=K\left(  i,j\right)  \text{.}%
\]
Hence $\mathbf{M}=\mathbf{K}$ and by (\ref{Kneadingdet}) and (\ref{FP15}) we
finally arrive at%
\[
\Delta=\det(\mathbf{I}-z\mathbf{K})=\det(\mathbf{I}-z\mathbf{M})=\det
(\mathbf{I}-z\mathbf{F}).
\]
This last relation is precisely what is stated in Theorem \ref{T1}.
\end{pf}

We now prove Theorem \ref{T2}.

\begin{pf}
[Proof of Theorem \ref{T2}]The idea is to regard a $Fib_{p}$ recurrence,
determined by a sequence of matrices $(\mathbf{A}_{n})_{n\in\mathbb{N}}$, as a
pair $(\varphi,\psi)$ of linear endomorphisms on the infinite dimensional
vector space $\boldsymbol{U}$ defined in (\ref{spaces}). This pair
$(\varphi,\psi)\in\mathrm{L}(\boldsymbol{U})\times\mathrm{L}(\boldsymbol{U})$
is now defined as follows:
\[
\varphi(\mathrm{u}_{0},\mathrm{u}_{1},\mathrm{u}_{2},...)=(\mathrm{0}%
,\mathrm{u}_{0},\mathrm{u}_{1},\mathrm{u}_{2},...)\text{ for all }%
(\mathrm{u}_{n})_{n\in\mathbb{N}}\in\boldsymbol{U},
\]
where $\mathrm{0}$ denotes the zero vector of $\mathbb{C}^{p}$, and%
\[
\psi(\mathrm{u}_{0},\mathrm{u}_{1},\mathrm{u}_{2},...)=(\mathrm{w}%
,\mathrm{u}_{0},\mathrm{u}_{1},\mathrm{u}_{2},...)\text{ for all }%
(\mathrm{u}_{n})_{n\in\mathbb{N}}\in\boldsymbol{U}\text{,}%
\]
with%
\[
\mathrm{w}=\sum\limits_{n\geq0}\mathbf{A}_{n}\mathrm{u}_{n}\in\mathbb{C}^{p}.
\]

Clearly, $(\varphi,\psi)$ is a pair of finite rank. Let us begin by proving
that $\Delta(\varphi,\psi)$ is actually the kneading determinant $\Delta$ of
the linear recurrence as defined in (\ref{Kneadingdet}).

Indeed, from the definitions of $\varphi$ and $\psi$ one has%
\[
\psi-\varphi=\omega_{1}\otimes\boldsymbol{e}_{1}+\omega_{2}\otimes
\boldsymbol{e}_{2}+\cdots+\omega_{p}\otimes\boldsymbol{e}_{p}\text{,}%
\]
where $\boldsymbol{e}_{i}\in\boldsymbol{U}$ denotes $i$-th vector of the
standard basis of $\boldsymbol{U}$ and $\omega_{i}\in\boldsymbol{U}^{\ast}$ is
the linear form defined by%
\[
\omega_{i}(\mathrm{u}_{0},\mathrm{u}_{1},\mathrm{u}_{2},...)=\sum
\limits_{n\geq0}\mathbf{R}_{n}(i)\mathrm{u}_{n},
\]
where $\mathbf{R}_{n}(i)=\left(  A_{n}(i,1)\cdots A_{n}(i,p)\right)
\in\mathbb{C}^{1\times p}$ denotes the $i$-th row of $\mathbf{A}_{n}$. By
Lemma \ref{LP1} we have%
\begin{equation}
\Delta(\varphi,\psi)=\det(\mathbf{I}-z\mathbf{M}), \label{FP16}%
\end{equation}
with%
\[
\mathbf{M}=\left(
\begin{array}
[c]{ccc}%
\underset{n\geq0}{\sum}\omega_{1}\varphi^{n}\left(  \boldsymbol{e}_{1}\right)
z^{n} & \cdots & \underset{n\geq0}{\sum}\omega_{1}\varphi^{n}\left(
\boldsymbol{e}_{p}\right)  z^{n}\\
\vdots & \ddots & \vdots\\
\underset{n\geq0}{\sum}\omega_{p}\varphi^{n}\left(  \boldsymbol{e}_{1}\right)
z^{n} & \cdots & \underset{n\geq0}{\sum}\omega_{p}\varphi^{n}\left(
\boldsymbol{e}_{p}\right)  z^{n}%
\end{array}
\right)  .
\]
But by (\ref{KneadMatrixK}), the $(i,j)$ entry of $\mathbf{M}$ is%
\[
\underset{n\geq0}{\sum}\omega_{i}\varphi^{n}\left(  \boldsymbol{e}_{j}\right)
z^{n}=\underset{n\geq0}{\sum}A_{n}(i,j)z^{n}=K\left(  i,j\right)  \text{.}%
\]
Hence $\mathbf{K}=\mathbf{M}$ and by (\ref{Kneadingdet}) and (\ref{FP16}) we
arrive at%
\begin{equation}
\Delta=\det(\mathbf{I}-z\mathbf{K})=\det(\mathbf{I}-z\mathbf{M})\text{.}
\label{FP11}%
\end{equation}
This formula is the first step in the proof of Theorem \ref{T2}. The second
step deals with the generating functions $G(\boldsymbol{u})$ of a $Fib_{p}$ recurrence.

Let $\pi:\boldsymbol{U}\rightarrow\mathbb{C}^{p}$ be the projection defined by
$\pi(\mathrm{u}_{0},\mathrm{u}_{1},\mathrm{u}_{2},...)=\mathrm{u}_{0}$. For
each $\alpha=1,...,p$, define the linear form $\pi_{\alpha}\in\boldsymbol{U}%
^{\ast}$, where $\pi_{\alpha}(\boldsymbol{u})$ is the $\alpha$-th coordinate
of $\pi(\boldsymbol{u})$ with respect the standard basis of $\mathbb{C}^{p}$.

Now let $\Theta(\boldsymbol{u})=(\mathrm{v}_{n})_{n\in\mathbb{Z}}$ be the
solution of the $Fib_{p}$ recurrence for the initial condition $\boldsymbol{u}%
=(\mathrm{u}_{n})_{n\in\mathbb{N}}\in\boldsymbol{U}$. Observe that from the
definition of $\psi$ one has%
\[
\pi\psi^{n}(\boldsymbol{u})=\mathrm{v}_{n}\text{ for }n\in\mathbb{N}\text{.}%
\]
Thus, the equalities%
\[
G(\boldsymbol{u})=\sum_{n\geq0}\mathrm{v}_{n}z^{n}=\sum_{n\geq0}\pi\psi
^{n}(\boldsymbol{u})z^{n}%
\]
and%
\begin{equation}
G_{\alpha}(\boldsymbol{u})=\sum_{n\geq0}\mathrm{v}_{n}^{(\alpha)}z^{n}%
=\sum_{n\geq0}\pi_{\alpha}\psi^{n}(\boldsymbol{u})z^{n}, \label{FP10}%
\end{equation}
hold for all $\boldsymbol{u}\in\boldsymbol{U}$ and $\alpha=1,...,p$.

Finally, we have all the elements to conclude the proof of Theorem \ref{T2}.

Let $\boldsymbol{e}_{\beta}$ be a vector of the standard basis of
$\boldsymbol{U}$ and $\alpha=1,...,p$. From Lemma \ref{LP3} and equality
(\ref{FP11}) one has%
\begin{align*}
z\sum_{n\geq0}\pi_{\alpha}\psi^{n}(\boldsymbol{e}_{\beta})z^{n}  &
=1-\frac{\det(\mathbf{I}-z\mathbf{M}_{\pi_{\alpha}}(\boldsymbol{e}_{\beta}%
))}{\det(\mathbf{I}-z\mathbf{M})}\\
&  =1-\Delta^{-1}\det(\mathbf{I}-z\mathbf{M}_{\pi_{\alpha}}(\boldsymbol{e}%
_{\beta})),
\end{align*}
with%
\[
\mathbf{M}_{\pi_{\alpha}}(\boldsymbol{e}_{\beta})=\left(
\begin{array}
[c]{cccc}%
\underset{n\geq0}{\sum}\omega_{1}\varphi^{n}\left(  \boldsymbol{e}_{1}\right)
z^{n} & \cdots & \underset{n\geq0}{\sum}\omega_{1}\varphi^{n}\left(
\boldsymbol{e}_{p}\right)  z^{n} & \underset{n\geq0}{\sum}\omega_{1}%
\varphi^{n}\left(  \boldsymbol{e}_{\beta}\right)  z^{n}\\
\vdots & \ddots & \vdots & \vdots\\
\underset{n\geq0}{\sum}\omega_{p}\varphi^{n}\left(  \boldsymbol{e}_{1}\right)
z^{n} & \cdots & \underset{n\geq0}{\sum}\omega_{p}\varphi^{n}\left(
\boldsymbol{e}_{p}\right)  z^{n} & \underset{n\geq0}{\sum}\omega_{p}%
\varphi^{n}\left(  \boldsymbol{e}_{\beta}\right)  z^{n}\\
\underset{n\geq0}{\sum}\pi_{\alpha}\varphi^{n}\left(  \boldsymbol{e}%
_{1}\right)  z^{n} & \cdots & \underset{n\geq0}{\sum}\pi_{\alpha}\varphi
^{n}\left(  \boldsymbol{e}_{p}\right)  z^{n} & \underset{n\geq0}{\sum}%
\pi_{\alpha}\varphi^{n}\left(  \boldsymbol{e}_{\beta}\right)  z^{n}%
\end{array}
\right)  .
\]
On the other hand, it is easy to see that $\mathbf{M}_{\pi_{\alpha}%
}(\boldsymbol{e}_{\beta})$ is actually the extended kneading matrix
$\mathbf{K}_{\alpha}(\beta)$ defined in (\ref{KneadMatrixKalfa}), hence%
\begin{align*}
z\sum_{n\geq0}\pi_{\alpha}\psi^{n}(\boldsymbol{e}_{\beta})z^{n}  &
=1-\Delta^{-1}\det(\mathbf{I}-z\mathbf{M}_{\pi_{\alpha}}(\boldsymbol{e}%
_{\beta}))\\
&  =1-\Delta^{-1}\det(\mathbf{I}-z\mathbf{K}_{\alpha}(\beta))\\
&  =1-\Delta^{-1}\Delta_{\alpha}(\beta)
\end{align*}
and by (\ref{FP10}) we finally arrive at%
\[
zG_{\alpha}(\boldsymbol{e}_{\beta})=1-\Delta^{-1}\Delta_{\alpha}(\beta).
\]
This last relation is precisely what is stated in Theorem \ref{T2}.
\end{pf}

Finally, we prove Theorem \ref{T3}.

\begin{pf}
[Proof of Theorem \ref{T3}]Let $\lambda_{1},...,\lambda_{m}$ be the dominant
eigenvalues of a $Fib_{p}$-recurrence whose kneading matrix $\mathbf{K}$ is
holomorphic on some open disk $D_{\rho}$ with $\rho>1.$ Therefore, the
kneading determinant $\Delta$ is holomorphic on $D_{\rho}$ and the zeros of
$\Delta$ lying in $\left\{  z\in\mathbb{C}:\left\vert z\right\vert
\leq1\right\}  $ are $z_{i}=\lambda_{i}^{-1}$, $i=1,...,m$.

From (\ref{KneadMatrixKalfa}), it is easy to see that every extended kneading
matrix $\mathbf{K}_{\alpha}\mathbf{(}\beta\mathbf{)}$ is also holomorphic on
$D_{\rho}$. Thus, every extended kneading determinant $\Delta_{\alpha}%
(\beta)=\det(\mathbf{I}-z\mathbf{K}_{\alpha}\mathbf{(}\beta\mathbf{)})$ is
holomorphic on $D_{\rho}$. By Theorem \ref{T2}, it turns clear that every
generating function%
\[
G_{\alpha}(\boldsymbol{e}_{\beta})=\frac{\Delta-\Delta_{\alpha}(\beta
)}{z\Delta}%
\]
is meromorphic on $D_{\rho}$. Moreover, as
\[
\Delta(0)-\Delta_{\alpha}(\beta)(0)=\det(\mathbf{I}-0\mathbf{K})-\det
(\mathbf{I}-0\mathbf{K}_{\alpha}\mathbf{(}\beta\mathbf{)})=1-1=0\text{,}%
\]
the meromorphic function $G_{\alpha}(\boldsymbol{e}_{\beta})$ has a removable
singularity at $0$. Hence, every pole of $G_{\alpha}(\boldsymbol{e}_{\beta})$
is a zero of $\Delta$.

This proves that the possible poles of $G_{\alpha}(\boldsymbol{e}_{\beta})$
lying in $\left\{  z\in\mathbb{C}:\left\vert z\right\vert \leq1\right\}  $ are
$z_{i}=\lambda_{i}^{-1}$,$i=1,...,m$.

Let us consider the Laurent's series of $G_{\alpha}(\boldsymbol{e}_{\beta})$
with respect to $z_{i}$%
\[
\sum_{j=-\mathrm{mul}(\lambda_{i})}^{+\infty}L(i,j)(z-z_{i})^{j}\text{.}%
\]
As each $z_{i}$ is a zero of $\Delta$ of order \textrm{mul}$(\lambda_{i})$,
the auxiliary function%
\begin{equation}
h(z)=G_{\alpha}(\boldsymbol{e}_{\beta})-\sum_{i=1}^{m}\sum_{j=-\mathrm{mul}%
(\lambda_{i})}^{-1}L(i,j)(z-z_{i})^{j} \label{FP25}%
\end{equation}
is holomorphic on some $D_{\rho^{\prime}}$, with $\rho>\rho^{\prime}>1.$
Consequently, the radius of convergence of $h(z)=%
{\textstyle\sum_{n\geq0}}
h_{n}z^{n}$ is grater than $1$ and one has
\begin{equation}
h_{n}\rightarrow0\text{.} \label{FP26}%
\end{equation}
On the other hand, combining (\ref{FP25}) with the formulas%
\[
(z-z_{i})^{j}=(z-\lambda_{i}^{-1})^{j}=\left(  -\lambda_{i}\right)  ^{-j}%
{\textstyle\sum_{n\geq0}}
\frac{(n-j-1)!}{\left(  -j-1\right)  !n!}\lambda_{i}^{n}z^{n},\text{ for }%
j\in\mathbb{Z}^{-},
\]
and%
\[
G_{\alpha}(\boldsymbol{u})=\sum\limits_{n\geq0}\mathrm{v}_{n}^{\left(
\alpha\right)  }z^{n}\text{,}%
\]
one gets%
\[
\mathrm{v}_{n}^{\left(  \alpha\right)  }-\sum_{i=1}^{m}\left(  \sum
_{j=1}^{\mathrm{mul}(\lambda_{i})}L(i,-j)\left(  -\lambda_{i}\right)
^{j}\frac{(n+j-1)!}{\left(  j-1\right)  !n!}\right)  \lambda_{i}^{n}%
=h_{n}\text{.}%
\]
Combining this last equality with (\ref{FP26}) and defining $c_{i,j}^{\left(
\alpha\right)  }(\boldsymbol{u})=L(i,-j)\left(  -\lambda_{i}\right)  ^{j}$, we
finally arrive at%
\[
\underset{n\rightarrow+\infty}{\lim}\left(  \mathrm{v}_{n}^{\left(
\alpha\right)  }-\sum_{i=1}^{m}\sum_{j=1}^{\mathrm{mul}(\lambda_{i})}%
\frac{c_{i,j}^{\left(  \alpha\right)  }(\boldsymbol{u})(n+j-1)!}{\left(
j-1\right)  !n!}\lambda_{i}^{n}\right)  =0\text{.}%
\]
This last relation is precisely what is stated in Theorem \ref{T3}.
\end{pf}

\textbf{Acknowledgement} We thank the valuable comments and suggestions from
the referee which improved the final version of this article. Partially funded
by FCT/Portugal through project PEst-OE/EEI/LA0009/2013 for CMAGDS.

\end{document}